\documentclass[11pt]{article}
\usepackage{amsthm,amssymb,amsmath,amscd
}
\usepackage[all]{xy}
\usepackage{mathrsfs}
\usepackage{multirow}
\usepackage{color}
\usepackage{url}

\usepackage{caption}
\usepackage{subfigure}
\usepackage{graphicx}
\usepackage{tikz-cd}\usepackage{enumitem}
\DeclareMathAlphabet{\mathcal}{OMS}{cmsy}{m}{n}
\DeclareSymbolFont{largesymbols}{OMX}{cmex}{m}{n}

\newtheorem{Def}{Definition}[section]
\newtheorem{Prop}[Def]{Proposition}
\newtheorem{Theo}[Def]{Theorem}
\newtheorem{Lem}[Def]{Lemma}
\newtheorem{Koro}[Def]{Corollary}

\newtheorem{Remark}[Def]{Remark}
\newtheorem{Example}[Def]{Example}

\newtheorem{Theox}{Theorem}

\newcommand{\indt}[3]{\begin{smallmatrix}#1\\#2\\#3
\end{smallmatrix}
}

\newcommand{\add}{\operatorname{add}\,}

\newcommand{\soc}{{\rm soc}}

\renewcommand{\top}{\operatorname{top}}

\newcommand{\opp}{^{\rm op}}

\newcommand{\findim}{\operatorname{findim}}
\newcommand{\pd}{{\rm pd }}
\newcommand{\id}{{\rm id}}

\newcommand{\Hom}{\operatorname{Hom}}
\newcommand{\rHom}{{\rm\bf R}{\rm Hom}}

\newcommand{\otimesL}{\otimes^{\rm\bf L}}

\newcommand{\End}{\operatorname{End}}
\newcommand{\Ext}{\operatorname{Ext}}
\newcommand{\Tor}{\operatorname{Tor}}

\newcommand{\D}[1]{\mathscr{D}(#1)}
\newcommand{\Dz}[1]{{\mathscr{D}}^+(#1)}

\newcommand{\Db}[1]{ \mathscr{D}^{\rm b}(#1)}

\newcommand{\Kb}[1]{ \mathscr{K}^{\rm b}(#1)}
\newcommand{\modcat}[1]{#1\mbox{{\rm -mod}}}

\newcommand{\Modcat}[1]{#1\mbox{{\rm -Mod}}}
\newcommand{\rModcat}[1]{\mbox{{\rm Mod-}}#1}

\newcommand{\pmodcat}[1]{#1\mbox{{\rm -proj}}}

\newcommand{\imodcat}[1]{#1\mbox{{\rm -inj}}}

\newcommand{\lra}{\longrightarrow}

\newcommand{\ra}{\rightarrow}

\def\scr#1{\mathscr{#1}}

\title{Stratified algebras, stratifying recollements and Iwanaga--Gorenstein algebras}
\author{Yiping Chen, Tiago Cruz, Steffen Koenig}
\date{}

\begin{document}
\maketitle

\renewcommand{\thefootnote}{\alph{footnote}}
\setcounter{footnote}{-1} 
\renewcommand{\thefootnote}{\alph{footnote}}
\setcounter{footnote}{-1} 

\begin{abstract}
  The standardly stratified algebras with finite global dimension are known to be exactly the quasi-hereditary algebras. In this article, we investigate the more general class of standardly stratified algebras that are Iwanaga--Gorenstein.
  For Iwanaga--Gorenstein algebras we prove that the concepts of properly stratified and standardly stratified coincide. Our main strategy is to characterise standardly stratified algebras via iterated half recollements of bounded derived categories and to prove that stratifying recollements of unbounded derived categories that restrict to half recollements of bounded derived categories can be lifted to full recollements of bounded derived categories when all the algebras involved are Iwanaga--Gorenstein.
 \end{abstract}


\section{Introduction}

In algebraic Lie theory and representation theory, highest weight categories
are described by quasi-hereditary algebras. Various more general
categories are described by
several kinds of standardly stratified algebras, with respect to a fixed
(partial or total) order of primitive idempotents.
Quasi-hereditary algebras are properly
stratified, and properly stratified algebras are both left standardly
stratified and right standardly stratified. In general, these classes of
algebras are quite different and the structural differences
between them have not yet been fully clarified. Moreover,
it is not clear under what assumptions some classes of stratified
algebras coincide. 

Quasi-hereditary algebras have finite global dimension, which distinguishes
them from non-quasi-hereditary stratified algebras. Another
widely used property of quasi-hereditary algebras is the existence of a
characteristic tilting module affording a double centraliser property called
Ringel duality; Ringel self-duality is a strong symmetry property. Classical references are \cite{ParshallScott1987, CPS1988b, DR1992, zbMATH04089744, Rin91}. 
For properly stratified algebras of infinite global dimension there exist a
characteristic tilting module and a characteristic cotilting module
and Ringel duality works smoothly when the characteristic tilting module
coincides with the characteristic cotilting module (see \cite{Agoston2000,
  FM2006}). This condition also plays a key role in the main theorem of
\cite{MO2004}
stating that the finitistic dimension of a properly
stratified algebra with a simple preserving duality is an even number,
when characteristic tilting and cotilting module coincide.
For properly stratified
algebras the two characteristic modules coincide
if and only if the algebra is Iwanaga--Gorenstein, that
is, its left and right regular module have finite injective dimension
(see \cite{CM2021}). Lie theoretic examples of properly stratified
algebras that are Gorenstein and whose characteristic tilting and
cotilting modules coincide include the blocks of {$\scr S$}-subcategories
in {$\scr O$}, see \cite[Sections 5 and 6]{FKM0} and \cite[4.5]{MO2004}.

Within the class of Iwanaga--Gorenstein algebras, the main result of
this article identifies the above three classes of stratified algebras
with each other.

\begin{Theox}[\bf Comparison Theorem \ref{mainresult}] \label{Theorem A}
	Let $R$ be a finite dimensional Iwanaga--Gorenstein
	algebra over a field, with a fixed ordering $e_1 < \ldots < e_n$
	on a complete set of pairwise orthogonal primitive idempotents. 
	
	Then $R$ is properly stratified with respect to the given ordering if and only
	if it is left standardly stratified with respect to the given ordering if
	and only if it is right standardly stratified with respect to the given
	ordering.
	
	In other words, within the class of Iwanaga--Gorenstein algebras the
	concepts of properly stratified, left standardly stratified and right
	standardly stratified coincide.
\end{Theox}

Like results in [CM21],
this result suggests that Iwanaga--Gorenstein algebras are a 
convenient generalisation of algebras of finite global dimension when studying
stratified algebras instead of quasi-hereditary algebras.

The main ingredient of the proof of the Comparison Theorem
uses another fundamental
property of stratified algebras, stratifications of derived categories by
iterated recollements. This property, already emphasised in
\cite{CPS1988b, ParshallScott1987}, seems to have been a
strong motivation for introducing quasi-hereditary algebras - which admit
iterated recollements by derived categories of vector spaces. Properly stratified structures and algebras, introduced by Dlab in \cite{Dl1}, likewise give rise to such stratifications via iterated recollements of bounded derived categories involving local algebras. This property
has, however, not been used frequently later on. Obstructions to the use of
recollements are that already the class of quasi-hereditary algebras is not
closed under derived equivalences (see for example \cite{LiuYang2016}) and that quasi-hereditary algebras sometimes
have other stratifications involving rather different recollements with
different factors. These further recollements are not stratifying and it
has been shown in \cite{AKLY2017a} that recollements are not in general
stratifying recollements. 

The proof of the Comparison Theorem illustrates that using particular kinds
of recollements instead of general ones is useful when working
with stratified algebras, which by
definition have stratifying recollements. Restricting attention to
stratifying recollements and even to half recollements that only use four of
the six functors allows to characterise the different kinds of stratified
algebras in the following way:

{\begin{Theox} [\bf Characterisation Theorem \ref{equi} and Corollary \ref{equiright}] \label{Theorem B}
		Let $R$ be a finite dimensional algebra over a
		field, with a fixed ordering $e_1 < \ldots <  e_n$ on a complete set of
		pairwise orthogonal primitive idempotents. Then:
		
		(a) $R$ is left standardly stratified with respect to the given ordering
		if and only if $R/Re_nR$ is left standardly stratified with respect to the
		induced ordering and there is an induced lower half stratifying  recollement
		involving the algebra $R$, a quotient  $R/Re_nR$ and a corner subalgebra $e_nRe_n$:
		
		\begin{equation} \label{d}
			\xymatrix@=3.5em{
				\Db{{R/Re_nR}}\ar[rr]^{R/Re_nR\otimesL_{R/Re_nR}-}&&\Db{{R}}\ar[rr]^{\rHom_R(Re_n, -)}\ar@/^2pc/[ll]^{\rHom_R(R/Re_nR, -)}&&\Db{{e_nRe_n}}\ar@/^2pc/[ll]^{\rHom_{e_nRe_n}(e_nR, -)}.}\end{equation}
		
		(b) $R$ is right standardly stratified with respect to the given ordering
		if and only if $R/Re_nR$ is right standardly stratified with respect to the
		induced ordering and there is an induced upper half stratifying  recollement:
		\begin{equation} 
			\xymatrix@=3.5em{
				\Db{{R/Re_nR}}\ar[rr]_{R/Re_nR\otimesL_{R/Re_nR}-}&&\Db{{R}}\ar[rr]_{\rHom_R(Re_n, -)}\ar@/_2pc/[ll]_{R/Re_nR\otimesL_R-}&&\Db{{e_nRe_n}}\ar@/_2pc/[ll]_{Re_n\otimesL_{e_nRe_n}-}.}\end{equation}
		
		(c) $R$ is properly stratified with respect to the given ordering if and
		only if $R/Re_nR$ is properly stratified with respect to the induced
		ordering and there is an induced (full) stratifying recollement:
		\begin{equation}
			\xymatrix@C=3.5em{
				\Db{\Modcat{R/Re_nR}}\ar[rr]^{R/Re_nR\otimesL_{R/Re_nR}-}&&\Db{\Modcat{R}}\ar@/^2pc/[ll]^{\rHom_R(R/Re_nR, -)}\ar@/_2pc/[ll]_{R/Re_nR\otimesL_R-}\ar[rr]^{\rHom_R(Re_n, -)}&&\Db{\Modcat{e_nRe_n}}\ar@/^2pc/[ll]^{\rHom_{e_nRe_n}(e_nR, -)}\ar@/_2pc/[ll]_{Re_n\otimesL_{e_nRe_n}-}.}\label{eq22}\end{equation}
		
	\end{Theox}
	
	Combining these characterisations with an extension of a previously obtained
	result from \cite{CK2017} on
        recollements involving Iwanaga--Gorenstein algebras we will be able to prove Theorem A.
	
	\bigskip
	
	This article is organised as follows. 
	In Section 2, we review known facts on standardly stratified algebras, introduce the notation to be used throughout the article, and present the definition of half recollements. Section 3 is divided into four parts: in the first we recall known facts about stratifying ideals; in the second we recall the concept of stratifying recollement in the sense of \cite{AKL2011}; in the third we discuss when the functors involved in these recollements restrict to bounded derived categories; and in the fourth we combine these results to prove Theorem B. In Section 4, we prove Theorem A. To this end, we first clarify how Iwanaga--Gorensteinness manifests itself in a recollement of bounded derived categories and we show that half recollements lift to full recollements when all algebras involved are Iwanaga--Gorenstein.
        We then apply these results to standardly stratified algebras in order to obtain Comparison Theorem A. Finally, in Section 5, we present examples, including a bounded derived category of a non-Iwanaga--Gorenstein algebra, which is stratified into bounded derived categories of Iwanaga--Gorenstein algebras via iterated recollements.
	
\section{Stratified algebras and recollements}

\paragraph{General conventions}
Let $k$ be a field. Assume that $R$ is a finite dimensional $k$-algebra.
 $\Modcat{R}$ is the category of left $R$-modules and $\modcat{R}$ its subcategory whose objects are finitely generated modules. 
 To distinguish left $R$-modules and right $R$-modules, we write $_{R}M$ and $M_R$, respectively. Given an $R$-module, we write $\add M$ (or $\add_R M$) to denote its additive closure. $D=\Hom_k(-, k)$ is the standard duality between the categories of finitely generated left and right  $R$-modules, while $R^{\opp}$ denotes the opposite algebra of $R$. The projective dimension of $M$ is denoted by $\pd_R M$ while the injective dimension is denoted by $\id_R M$. By $\pmodcat{R}$ we denote the subcategory consisting of finitely generated projective $R$-modules and by $\imodcat{R}$ the subcategory of $\modcat{R}$ consisting of finitely generated injective modules.
 Given a finite set $\Theta$ of $R$-modules, $\mathcal{F}(\Theta)$ is the full subcategory of $\modcat{R}$ whose objects admit a finite filtration whose factors belong to $\Theta$. 
We read the arrows in a path algebra like morphisms, that is, from right to left.

A complex of left $R$-modules $X^\bullet$ is a sequence of left $R$-modules and left $R$-module homomorphisms $d^i\colon X^i\rightarrow X^{i+1}$ such that $d^{i+1}\circ d^i=0$. A complex $X^\bullet$ is called bounded below (resp. bounded above) if $X^i=0$ for all but finitely many $i<0$ (resp. $i>0$), and bounded if it is bounded below and bounded above.
By $\D{R}$ (resp. $\Dz{R}$) we denote the derived category of complexes (resp. bounded below complexes) of $R$-modules, and by $\Db{R}$ the bounded derived category of $R$, that is, the derived category of bounded complexes of $R$-modules.
We write $\Kb{\pmodcat{R}}$ (resp. $\Kb{\imodcat{R}}$) for the bounded homotopy category
of finitely generated projective (resp. injective) $R$-modules. All these five  categories are triangulated, with shift functor $[1]$ given by
$(X[1])^n:=X^{n+1}$.
As usual, the derived functors of $\Hom_R(-, -)$ and $-\otimes_R -$ will be denoted by $\rHom_R(-, -)$ and $-\otimesL_R -$, respectively.

\subsection{Standardly and properly stratified algebras}

Let $R$ be a finite dimensional $k$-algebra and ${\bf e }:=\{e_1, \ldots, e_n\}$ a complete set of primitive pairwise orthogonal idempotents. Without loss of generality, we may assume that these idempotents are ordered according to the linear order $1<2<\ldots<n$. Hence, $P(i)=Re_i$ represent the indecomposable projective $R$-modules and $S(i)=\top Re_i$ represent the simple $R$-modules for $i=1, \ldots, n$, up to isomorphism.   Recall that the (left) standard module $\Delta(i)$ is the maximal quotient of $P(i)$ without composition factors $S(j)$ with $j>i$. Hence, $\Delta(i)=P(i)/R(e_{i+1}+\cdots+ e_n)Re_i$. Analogously, the (right) standard module $\Delta_{\opp}(i)$ is the module $e_iR/e_iR(e_{i+1}+\cdots+e_n)R$.
By $\Delta$ we denote the set $\{ \Delta(i): 1 \leq i \leq n \}$ of
standard modules; $\Delta_{\opp}$ is defined analogously.

\begin{Def}
	The pair $(R, {\bf e}=\{e_1, \ldots, e_n\})$ is called 
	\begin{itemize}
		\item a left standardly stratified algebra if ${}_R R\in \mathcal{F}(\Delta)$;
		\item a right standardly stratified algebra if $R_R\in \mathcal{F}(\Delta_{\opp})$;
		\item a properly stratified algebra if it is right and left standardly stratified.
	\end{itemize}
\end{Def} 

As the name suggests the concept of standardly stratified algebras is not left-right symmetric.
An explicit example of a left standardly stratified algebra which is not properly stratified is the following:

\vspace{-0.5cm}

\begin{Example}
	Let $R$ be the quotient of the path algebra of the quiver $
	\begin{tikzcd}
		1  & 2 \arrow[l, "\alpha", swap] \arrow[loop, distance=3em, out=35, in=-35, "\beta"]
              \end{tikzcd}$ \vspace{-1em}
              modulo the relations $\beta^2=\alpha\beta=0$. Denote by $e_i$ the primitive idempotent associated with the vertex $i$ for $i=1, 2$. Then $(R, \{e_1, e_2\})$ is left standardly stratified since $\Delta(2)=Re_2$ and $\Delta(1)=Re_1$ with the ordering $1<2$. It is not, however, right standardly stratified since $e_1R$ and $e_2R$ have Loewy length two with $\soc(e_1R)=\top(e_2R)$ and so $e_1R$ cannot be obtained by extensions of $\Delta_{\opp}(1)=\top(e_1R)$ and $\Delta_{\opp}(2)=e_2R$.
\end{Example}

 The labelling matters because with another labelling the algebra $R$ of the previous example becomes properly stratified (see \cite[Example 2.6]{CM2021}).

\begin{Remark}
	\begin{enumerate}
        \item Some sources use a partial order to label the standard modules instead of a linear order as we do here. The two definitions are equivalent since every partial order can be extended to a linear order (see
\cite[Statements 1.2.3 and 1.2.4]{Cou2020}, \cite[Proposition 1.3.]{UY1990}).
		\item The terminology for standardly stratified algebras has varied in the literature. Left standardly stratified algebras were also called strongly standardly stratified algebras in \cite{FM2006} to distinguish them from the more general standardly stratified algebras studied in \cite{Cline1996}. The latter are now known in the literature as pre-standardly stratified algebras since for those the labelling set only needs to have a partial preorder (see for instance \cite{Mendoza2009}).
	\end{enumerate}
	
\end{Remark}

Properly stratified algebras include local algebras, hereditary algebras, quasi-hereditary algebras, centraliser algebras of nilpotent matrices and many algebras appearing in algebraic Lie theory. This includes, for instance, algebras related with certain subcategories of Harish-Chandra bimodules and generalisations of the category $\mathcal{O}$ of a complex semi-simple  Lie algebra (see for example \cite{FKM0, FKM1, FKM2, FM2006} and \cite{CM2021}). 

In the context of standardly stratified algebras, quasi-hereditary algebras appear as a distinguished subclass, for which the standard modules are already Schurian modules. In the general case, the Schurian counterparts of the standard modules are the proper standard modules.
The (left) proper standard module $\overline{\Delta}(i)$ is then the maximal quotient module of $\Delta(i)$ so that the multiplicity of $S(i)$ in $\overline{\Delta}(i)$ is exactly one. Dually, there are the (right) proper standard modules $\overline{\Delta}_{\opp}(i)$.

Standardly stratified algebras also come equipped with further classes of modules: costandard, proper costandard and, depending on the side, special tilting or cotilting modules that respect the stratification. For each $i=1, \ldots, n$, we write $\nabla(i)=D\Delta_{\opp}(i)$ and $\overline{\nabla}(i)=D\overline{\Delta}_{\opp}(i)$ to denote the (left) costandard and (left) proper costandard indexed by $i$. By $\overline{\Delta}$ we denote the set of proper standard modules.
The notations $\overline{\nabla}$, $\overline{\Delta}_{\opp}$, and their variants are defined analogously.

The next theorem collects some interactions between between left and right (proper) standard modules and their properties.

\begin{Theo}  Let $R$ be a finite dimensional algebra equipped with a complete linearly ordered set of primitive pairwise orthogonal idempotents ${\bf e}$.
  Then:
	\begin{enumerate}
		\item The following assertions are equivalent.
		\begin{enumerate}
			\item $(R, {\bf e})$ is a left standardly stratified algebra;
			\item $R_R\in \mathcal{F}(\overline{\Delta}_{\opp})$;
			\item ${}_R DR\in \mathcal{F}(\overline{\nabla})$;
			\item $(DR)_R\in \mathcal{F}(\nabla_{\opp})$.
		\end{enumerate}
		\item If the above statements hold, then the following properties are satisfied:
		\begin{enumerate}
			\item $\Ext_R^{i>0}(\mathcal{F}(\Delta), \mathcal{F}(\overline{\nabla}))=0$.
			\item There exists a unique multiplicity-free tilting left $R$-module $T$ such that $\mathcal{F}(\overline{\nabla})=T^{\perp}$. In particular, $\add T=\mathcal{F}(\Delta)\cap \mathcal{F}(\overline{\nabla})$ and for each $i=1, \ldots, n$ there are exact sequences $$0\rightarrow \Delta(i)\rightarrow T(i)\rightarrow X(i)\rightarrow 0$$ with $X(i)\in \mathcal{F}(\Delta(1), \ldots, \Delta(i-1))$ and $T=\oplus_{i=1}^n T(i)$ for indecomposable $R$-modules $T(i)$.
			\item $R$ and $R^{\opp}$ have finite finitistic dimension.
		\end{enumerate}
	\end{enumerate}
\end{Theo}
\begin{proof}
	For the first statement see \cite[Theorem 1.1]{Agoston2000}, \cite{Lakatos} and \cite{ADL98}. 2(a) is straightforward, see for example \cite[Proposition 1.2]{Platzeck2001} or \cite[Lemma 1.2]{Agoston2000} and the references therein. For 2(b) see \cite[Theorem 2.1, Proposition 2.2, Lemma 2.5]{Agoston2000}. For 2(c) see \cite[Theorem 2.1, Theorem 3.1]{Agoston2000a}.
\end{proof}

This collection of results shows that, in practice, the existence of both left and right standardly stratified structures can be verified using left modules.
 Indeed it follows that $(R, {\bf e})$ is right standardly stratified if and only if ${}_R R\in \mathcal{F}(\overline{\Delta})$. A major difference between left and right standardly stratified algebras lies in the type of characteristic modules they admit: while a left standardly stratified algebra possesses a characteristic tilting module, a right standardly stratified algebra possesses a characteristic cotilting module.
 
 A module that satisfies 2(b) is called \emph{characteristic tilting module}. Dually, a cotilting module $C$ that satisfies $\mathcal{F}(\overline{\Delta})={}^\perp C$ is a \emph{characteristic cotilting module}. 
 Even though properly stratified algebras admit both a characteristic cotilting and a characteristic tilting module, the two characteristic modules need not coincide (see for instance \cite{FM2006}). In many examples coming from Lie theory they do coincide like in the case of quasi-hereditary algebras.

 When characteristic tilting and cotilting module do coincide, the algebras exhibit nice symmetries and homological properties and appear to be particularly well-behaved from the perspective of Ringel duality. In general, however, the endomorphism algebra of a characteristic tilting module over a properly stratified algebra is not necessarily properly stratified. On the positive side, if the characteristic tilting is the same as the characteristic cotilting module, then the Ringel dual of a properly stratified algebra is again properly stratified (see \cite[Theorem 5]{FM2006}). This condition also plays a crucial role in \cite{MO2004}, where computing the finitistic dimension of properly stratified algebras with a simple preserving duality  is reduced to determining the projective dimension of the characteristic tilting module.

\begin{Theo} (\cite[Theorem 2.2]{CM2021})  \label{maintheoremCM2021}
	Let $(R, {\bf e})$ be a properly stratified algebra with characteristic tilting module $T$. Then the following assertions are equivalent.
	\begin{enumerate}
		\item $R$ is Iwanaga--Gorenstein.
		\item $T$ is a characteristic cotilting module.
		\item The modules in $\mathcal{F}(\nabla)$ are exactly those that have a finite resolution by modules in $\add T$.
	\end{enumerate}
\end{Theo}

In \cite{zbMATH00871351} (see also \cite{Agoston2000}), it has been shown that
the standardly stratified algebras with finite global dimension are exactly
the quasi-hereditary algebras. Hence, \cite[Theorem 2.2]{CM2021} suggests to
look for similar results involving Gorensteinness.
To tackle such questions, we will use recollements of
triangulated categories as the major tool.

\subsection{Half recollements}

Recollements of triangulated categories provide a framework to deconstruct triangulated categories into smaller ones, similar to deconstructing a module occuring as middle term of a short exact sequence. For instance, every derived category of a quasi-hereditary algebra is obtained by glueing derived categories of semi-simple algebras via recollements. For standardly stratified algebras, we will see that it is useful to consider half recollements of triangulated categories.

Let $\mathcal{T}, \mathcal{T}_1, \mathcal{T}_2$ be triangulated categories. An upper half recollement of triangulated categories is a diagram
$$\xymatrix{
	\mathcal{T}_1\ar[rr]_{i_*}&&\mathcal{T}\ar@/_1.5pc/[ll]_{i^*}\ar[rr]_{j^!}&&\mathcal{T}_2\ar@/_1.5pc/[ll]_{j_!}}$$
of triangle functors such that
\begin{enumerate}
	\item $(i^*, i_*)$ and $(j_!, j^!)$ are adjoint pairs;
	\item $i_*$ and $j_!$ are full embeddings;
	\item $i^*\circ j_!=0$ and $j^!\circ i_*=0$;
	\item for each $X\in \mathcal{T}$ there exists a triangle
	$j_!j^!(X)\rightarrow X\rightarrow i_*i^*(X)\rightarrow j_!j^!(X)[1]. $
\end{enumerate}
Dually, a lower half recollement of triangulated categories is a diagram
$$\xymatrix{
	\mathcal{T}_1\ar[rr]^{i_!}&&\mathcal{T}\ar@/^1.5pc/[ll]^{i^!}\ar[rr]^{j^*}&&\mathcal{T}_2\ar@/^1.5pc/[ll]^{j_*}}$$ of triangle functors such that
\begin{enumerate}
	\item $(i_!, i^!)$ and $(j^*, j_*)$ are adjoint pairs;
	\item $j_*$ and $i_!$ are full embeddings;
	\item $i^!\circ j_*=0$ and $j^*\circ i_!=0$;
	\item for each $X\in \mathcal{T}$ there exists a triangle
	$i_!i^!(X)\rightarrow X\rightarrow j_*j^*(X)\rightarrow i_!i^!(X)[1]. $
\end{enumerate}
A recollement of triangulated categories in the sense of \cite{Beilinson1982}
is a diagram 
$$\xymatrix{
  \mathcal{T}_1\ar[rr]^{i_*=i_!}&&\mathcal{T}\ar@/^1.5pc/[ll]^{i^!}\ar@/_1.5pc/[ll]_{i^*}\ar[rr]^{j^!=j^*}&&\mathcal{T}_2\ar@/_1.5pc/[ll]_{j_!}\ar@/^1.5pc/[ll]^{j_*}}$$ such that the functors $(i^*, i_*, j_!, j^!)$ fit into an upper half recollement and the functors $(i_!, i^!, j^*, j_*)$ fit into a lower half recollement.

In the literature on abelian or triangulated categories,
lower or upper half recollements, respectively,
sometimes are called localisations or colocalisations, respectively,
of abelian or triangulated categories.

\section{Characterising standardly stratified algebras via half recollements}

Dlab in \cite{Dl1} introduced properly stratified algebras by defining them recursively in terms of primitive idempotents that induce recollements of bounded derived categories; the specific functors defining these recollements were not made explicit. 
In a similar vein, we show in this subsection that both left and right standardly stratified algebras can be characterised recursively using primitive idempotents that give rise to particular half recollements of bounded derived categories,
with specified functors.

Our approach recovers the original definition from \cite{Dl1} while making explicit that the recollements involved are indeed stratifying recollements (in the sense of \cite{AKLY2017a}) of bounded derived categories.

The aim of this section is to prove the following theorem.

\begin{Theo}\label{equi}
	Let $(R, \{e_1, e_2, \ldots, e_n\})$ be a finite dimensional $k$-algebra, and $\{e_i\mid 1 \leqslant i\leqslant n\}$ be a complete set of primitive pairwise orthogonal idempotents. Denote by $\overline{e}_i$ the image of $e_i$ in $R/Re_nR$. Then the following assertions are equivalent.
		\begin{enumerate}[label=(a\arabic*)]
			\item
			$(R, \{e_1, e_2, \ldots, e_n\})$ is a left standardly stratified algebra with $1 < \cdots < n$.
			\item
			$(R/Re_nR, \{\overline{e}_1, \overline{e}_2, \ldots, \overline{e}_{n-1}\})$ is a left standardly stratified algebra, $Re_nR$ is a projective left $R$-module. 
			\item
			$(R/Re_nR, \{\overline{e}_1, \overline{e}_2, \ldots, \overline{e}_{n-1}\})$ is a left standardly stratified algebra,
			and there is a lower half recollement
			\begin{equation} \label{d}
				\xymatrix@=3.5em{
					\Db{{R/Re_nR}}\ar[rr]^{R/Re_nR\otimesL_{R/Re_nR}-}&&\Db{{R}}\ar[rr]^{\rHom_R(Re_n, -)}\ar@/^2pc/[ll]^{\rHom_R(R/Re_nR, -)}&&\Db{{e_nRe_n}}\ar@/^2pc/[ll]^{\rHom_{e_nRe_n}(e_nR, -)}.}\end{equation}
			\item
			$(R/Re_nR, \{\overline{e}_1, \overline{e}_2, \ldots, \overline{e}_{n-1}\})$ is a left standardly stratified algebra,
			and there is an upper half recollement
			\begin{equation}\label{u}
				\xymatrix@=3.2em{
					\Db{(R/Re_nR) ^{\opp}}\ar[rr]_{-\otimesL_{R/Re_nR}R/Re_nR}&&\Db{{R^{\opp}}}\ar@/_2pc/[ll]_{-\otimesL_R R/Re_nR}\ar[rr]_{\rHom_{R^{\opp}}(e_nR, -)}&&\Db{(e_nRe_n)^{\opp}}\ar@/_2pc/[ll]_{-\otimesL_{e_nRe_n}e_nR}.}\end{equation}\end{enumerate}
\end{Theo}

\bigskip

The proof of Theorem \ref{equi} will be based on the following three key facts:

\begin{enumerate}
\item[(1)] An algebra $(R, \{e_1, \ldots, e_n\})$ is left standardly stratified if and only if $Re_nR$ is a left projective $R$-module and $R/Re_nR$ is left standardly stratified (see Subsection \ref{Stratifying ideals}).
	\item[(2)] A stratifying ideal $ReR$ induces a recollement relating the (unbounded) derived categories of $R$, $R/ReR$ and $eRe$ and such a recollement is precisely the derived version of the recollement relating the module categories of $R$, $R/ReR$ and $eRe$ (see Subsection \ref{Stratifying recollements of unbounded derived categories}).
	\item[(3)] A stratifying ideal $ReR$ induces a (lower) half recollement relating the bounded derived categories of $R$, $R/ReR$ and $eRe$ exactly when $eR$ has finite projective dimension as left $eRe$-module (see Subsection \ref{half recollements subsection}).
\end{enumerate} 

The fact (2) is nowadays considered folklore. One proof can be found in \cite[Example 4.5(4)]{AKL2011} using the theory of homological ring epimorphisms.
To establish (3), we combine (2) with the fact that a triangle functor between unbounded derived categories restricts to the bounded homotopy categories of projective modules whenever it admits a right adjoint that restricts to the bounded derived categories. Since local algebras have finitistic dimension zero, $eR$ is projective if and only if it has finite projective dimension as $eRe$-module when $e$ is a primitive idempotent. Thus, Theorem \ref{equi} will be proved in
Subsection \ref{subsection-characterisations} by combining the facts (1) and (3).

In particular, our approach shows that if the quotient algebra
$R/ReR$ is left standardly stratified such that $R\rightarrow R/ReR$ induces a half recollement of bounded derived categories, then $ReR$ is also a particular case of a stratifying ideal.

\subsection{Stratifying ideals} \label{Stratifying ideals}

Quasi-hereditary algebras continue to serve, to some extent, as guiding examples for identifying desirable properties of standardly stratified algebras. A prominent feature of quasi-hereditary algebras
is that they can be characterised from a ring-theoretic perspective by heredity chains. In the context of standardly stratified algebras, heredity ideals are replaced by particular stratifying ideals.  

While the term stratifying ideal was formally introduced in \cite{Cline1996}, the underlying concept has appeared also in \cite{Auslander1992b} under the name of strong idempotent ideal (cf.$\!$ (d) below). 

\begin{Def}\label{stratifyingidealdef}
	Let $e\neq 0$ be an idempotent of $R$. We say that $ReR$ is a stratifying ideal of $R$ (see \cite[(2.1)]{Parshall1989}, \cite[2.1.1, 2.1.2]{Cline1996}) if one of the following equivalent conditions holds:
	\begin{enumerate}[label=(\alph*)]
		\item
		The induced functor of derived categories $\iota_*\colon {}\Dz{\Modcat{R/ReR}}\ra \Dz{\Modcat{R}}$ is an embedding;
		\item
		$Re\otimesL_{eRe}eR\cong ReR$;
		\item
		Multiplication induces an isomorphism $Re\otimes_{eRe}eR\ra ReR$ of $R$-bimodules and ${\Tor^{eRe}_n(Re, eR)=0}$ for all $n\geqslant 1$;
		\item For all left $R/ReR$-modules $M, M'$, the natural map $\Ext_{R/ReR}^i(M, M')\rightarrow \Ext_R^i(M, M')$ is an isomorphism for all $i\geq 0$.
	\end{enumerate}
\end{Def}

The following well-known assertion
shows that an idempotent ideal $ReR$ is already stratifying
when it is projective as a left or right module.

\begin{Remark}\label{progen} Let $R$ be a finite dimensional $k$-algebra, and $e$ be a primitive idempotent of $R$. Suppose that one of the following conditions holds:
	\begin{enumerate}[label=(\alph*)]
		\item $ReR$ is a projective left $R$-module.
		\item $ReR$ is a projective right $R$-module.
	\end{enumerate}
	Then the ideal $ReR$ is a stratifying ideal of $R$.  
\end{Remark}
\begin{proof} See \cite[2.1.2(b)]{Cline1996}.
\end{proof}

\begin{Remark}
It is tempting to consider another two conditions that look related to
those in the previous remark: $eR$ is a projective left $eRe$-module, or
$Re$ is a projective right $eRe$-module. These conditions are, however, not
equivalent to those in the previous remark,
as the following example shows. Let $R$ be a
selfinjective Nakayama algebra with two simples and vanishing radical square,
over a field $k$. Let $e=e_1$ be a primitive idempotent. Then
$Re \otimes_k eR$ has dimension four over $k$, while the dimension of $ReR$
is three. Hence these two bimodules are not isomorphic and thus $ReR$ is not
a stratifying ideal. This example also shows that the multplication map
$Re \otimes_k eR \rightarrow ReR$ need not be injective.
\end{Remark}

Stratifying ideals appear naturally in the theory of standardly stratified algebras. This is a consequence of the following.

  \begin{Lem} \label{CM2021} The following statements hold.
	\begin{enumerate}[label=(\alph*)]
		\item If $(R, \{e_1, \ldots, e_n\})$ is a right standardly stratified algebra, then $Re_n$ is a projective generator of $e_nRe_n$ and $Re_nR$ is projective as a right $R$-module.
		\item If $(R, \{e_1, \ldots, e_n\})$ is a left standardly stratified algebra, then $e_nR$ is a projective generator of $e_nRe_n$ and $Re_nR$ is projective as a left $R$-module.
	\end{enumerate}
\end{Lem}
\begin{proof}
  For (a), see \cite[Lemma 1.6]{CM2021}. (b) follows from (a) by recalling
  that $(R, \{e_1, \cdots, e_n\})$ is left standardly stratified if and only
  if its opposite algebra is right standardly stratified.
\end{proof}

So by the previous two results it follows that if $(R, {\bf e})$ is a left or
a right standardly stratified algebra, then $Re_nR$ is a stratifying ideal. 
Moreover, in the left standardly stratified case $\add{Re_nR} =
\add{Re_n}$ and in the right standardly stratified case $\add{Re_nR_R} =
\add{e_nR}$.

\begin{Koro} \label{coro1}
	Let $ReR$ be a stratifying ideal of $R$.
	The following assertions hold.
	\begin{enumerate}[label=(\alph*)]
		\item $ReR$ is projective as a left $R$-module if and only if $eR$ is projective as a left $eRe$-module. Dually, $ReR$ is projective as a right $R$-module if and only if $Re$ is projective as a right $eRe$-module.  
		\item If $e$ is a primitive idempotent of $R$, then the projective dimensions of $ReR$ as a left $R$-module and a right $R$-module are $0$ or $\infty$.  
	\end{enumerate}
\end{Koro}
\begin{proof}
	For (a), we refer to \cite[Remark 3.2]{Koenig2009}. 
	
	Note that the finitistic dimension of a local algebra is zero (see for example \cite[Lemma 2.1]{CM2021}). Assume that $e$ is a primitive idempotent of $R$. Then the projective dimension of $eR$ as a left $eRe$-module is $0$ or $\infty$. By the derived tensor-hom adjunction, $\Hom_{\scriptsize \D{\Modcat{R}}}(ReR, S[i])\cong \Hom_{\scriptsize \D{\Modcat{R}}}(Re\otimesL_{eRe} eR, S[i])\cong \Hom_{\scriptsize \D{\Modcat{eRe}}}(eR, eS[i])$ for every simple $R$-module and $i\in \mathbb{N}$. Since every simple $eRe$-module is of the form $eS$ (see for instance \cite[Theorem 6.2g]{zbMATH05080041}), it then follows that the projective dimension of $ReR$ as a left $R$-module is $0$ or $\infty$. A similar argument holds for $ReR$ as a right $R$-module. 
\end{proof}

A consequence of the previous two results is that an algebra $(R, \{e_1, \ldots, e_n\})$ is left standardly stratified if and only if $Re_nR$ is a stratifying ideal which is a projective left $R$-module and $R/Re_nR$ is left standardly stratified with respect to the complete set of primitive idempotents given by the images of $\{e_1, \ldots, e_{n-1}\}$ in $R/Re_nR$.
This reflects the inductive nature of standardly stratified algebras, which
can be described recursively in terms of local algebras.

\subsection{Stratifying recollements of unbounded derived categories} \label{Stratifying recollements of unbounded derived categories}

Let $R$ be a finite dimensional algebra over a field. 
In \cite{zbMATH06346473}, Psaroudakis and Vit\'oria established
that every recollement of abelian categories with $\Modcat{R}$ in the middle is equivalent to a recollement of the form 
\begin{equation}
	\xymatrix@C=3.5em{
			\Modcat{R/ReR}\ar[rr]^{R/ReR\otimes_{R/ReR}-}&& \Modcat{R}\ar@/^2pc/[ll]^{\Hom_R(R/ReR, -)}\ar@/_2pc/[ll]_{R/ReR\otimes_R-}\ar[rr]^{\Hom_R(Re, -)}&&\Modcat{eRe}\ar@/^2pc/[ll]^{\Hom_{eRe}(eR, -)}\ar@/_2pc/[ll]_{Re\otimes_{eRe}-}}, \label{modulecatrecollement}
\end{equation}
where $e$ is an idempotent of $R$. In general, this recollement does not lift to the level of bounded derived categories or the level of unbounded derived categories, unless the idempotent involved satisfies additional properties. Indeed, Cline, Parshall and Scott showed that in the context of quasi-hereditary algebras if $ReR$ is a heredity ideal, then (\ref{modulecatrecollement}) induces a recollement of bounded derived categories (see for example \cite{ParshallScott1987} or \cite{CPS1988b} and also \cite{Koenig1991}). Later, Miyachi proved in \cite{Miyachi1992} that this occurs if and only if $ReR$ is a recollement ideal, a condition stronger than that of being stratifying. 
The importance of stratifying ideals lies in the fact that for those the recollement of abelian categories (\ref{modulecatrecollement}) lifts to one of unbounded derived categories.
	\begin{Lem}\label{main}Let $R$ be a finite dimensional $k$-algebra, and $e$ be an idempotent of $R$. Then $ReR$ is a stratifying ideal of $R$ if and only if there is the following recollement
	\begin{equation}\label{eq:7}
		\xymatrix@C=3.5em{
			\D{\Modcat{R/ReR}}\ar[rr]^{R/ReR\otimesL_{R/ReR}-}&&\D{\Modcat{R}}\ar@/^2pc/[ll]^{\rHom_R(R/ReR, -)}\ar@/_2pc/[ll]_{R/ReR\otimesL_R-}\ar[rr]^{\rHom_R(Re, -)}&&\D{\Modcat{eRe}}\ar@/^2pc/[ll]^{\rHom_{eRe}(eR, -)}\ar@/_2pc/[ll]_{Re\otimesL_{eRe}-}.}\end{equation}
\end{Lem}
\begin{proof}
	If such a recollement exists, then it is clear by Definition \ref{stratifyingidealdef} that $ReR$ is stratifying. For the other implication, see for example \cite[Example 4.5(4)]{AKL2011}.
\end{proof}
Therefore, recollements of the form (\ref{eq:7}) are known as {\em stratifying recollements}.
 
Condition (c) of Definition \ref{stratifyingidealdef} is left-right symmetric, that is $ReR$ is a stratifying ideal of $R$ if and only if $ReR=R^{\opp}eR^{\opp}$ is a stratifying ideal of $R^{\opp}$. Consequently, if $ReR$ is a stratifying ideal of $R$, then 
there exists also a recollement 
\begin{equation}\label{dual of 7}
	\xymatrix@C=3.5em{
		\D{\rModcat{R/ReR}}\ar[rr]^{-\otimesL_{R/ReR} R/ReR}&&\D{\rModcat{R}}\ar@/^2pc/[ll]^{\rHom_{R^{\opp}}(R/ReR, -)}\ar@/_2pc/[ll]_{-\otimesL_R R/ReR}\ar[rr]^{\rHom_{R^{\opp}}(eR, -)}&&\D{\rModcat{eRe}}\ar@/^2pc/[ll]^{\rHom_{(eRe)^{\opp}}(Re, -)}\ar@/_2pc/[ll]_{-\otimesL_{eRe}eR}.}\end{equation}

\subsection{Half stratifying recollements of bounded derived categories}\label{half recollements subsection}

In general, the stratifying recollements (\ref{eq:7}) and (\ref{dual of 7})
do not restrict to recollements involving bounded derived module categories $\Db{{\rm mod}}$ of finitely generated modules. In what follows, we investigate necessary and sufficient conditions under which a recollement restricts to $\Db{{\rm mod}}$.

 To this end, we begin by recalling how to detect whether triangle functors restrict to $\Db{{\rm mod}}$ via their adjoint functors.

\subsubsection{Restrictions to $\Db{\operatorname{mod}}$}

The following lemma is taken from \cite[Lemma 2.7]{AKLY2017} and \cite[Theorem 4.1 and Theorem 5.1]{Neeman1996}, see also \cite[Lemma 2.2]{CHQW2023}. 

\begin{Lem} \label{adjoint}Let $R_1$ and $R_2$ be finite dimensional $k$-algebras. Let $F:{}\D{\Modcat{R_1}}\ra \D{\Modcat{R_2}}$ be a triangle functor with a right adjoint $G$. Then the following statements are equivalent.
	
	\begin{enumerate}[label=(\alph*)]
		\item
		$F$ restricts to $\Kb{{\rm proj}}$.
		\item
		$G$ admits a right adjoint.
		\item
		$G$ restricts to $\Db{{\rm mod}}$.
	\end{enumerate}
\end{Lem}

The dual formulation of Lemma~\ref{adjoint} will be essential for our purposes.

\begin{Lem}\label{adjointdual}
	Let $R_1$ and $R_2$ be finite dimensional $k$-algebras. Let $H\colon \D{\Modcat{R_1}}\ra \D{\Modcat{R_2}}$ be a triangle functor with a right adjoint $J$. Then the following assertions are equivalent.
	\begin{enumerate}[label=(\alph*)]
		\item $H$ restricts to $\Db{{\rm mod}}$.
		\item $J$ restricts to $\Kb{{\rm inj}}$.
	\end{enumerate}
\end{Lem}
\begin{proof}
The argument is essentially the same as the one used for lemma \ref{adjoint},
but by characterising $D^b$ in terms of injectives instead of projectives,
see also \cite[Lemma 4]{YY2016}.
\end{proof}

The following provides an easy to verify criterion to determine when two
of the embeddings in a stratifying recollement restrict to $\Db{{\rm mod}}$.

\begin{Lem} \label{lemmanew2}Let $R$ be a finite dimensional algebra, and $e$ be an idempotent. Then the following hold.
	\begin{enumerate}[label=(\alph*)]
		\item 
		The functor $-\otimesL_{eRe} eR:{} \D{\rModcat{eRe}}\ra \D{\rModcat{R}}$ restricts to $\Db{\rm mod}$ if and only if $eR$ has finite projective dimension as an $eRe$-module.
		\item 
		The functor $\rHom_{eRe}(eR, -):{} \D{\Modcat{eRe}}\ra \D{\Modcat{R}}$ restricts to $\Db{\rm mod}$ if and only if $eR$ has finite projective dimension as an $eRe$-module.
	\end{enumerate}
\end{Lem}
\begin{proof} It is clear that $-\otimesL_{eRe} eR:{} \D{\rModcat{eRe^{\opp}}}\ra \D{\rModcat{R^{\opp}}}$ restricts to $\Db{\rm mod}$ when $eR$ has finite projective dimension as a left $eRe$-module.
	Assume that $-\otimesL_{eRe} eR$ restricts to $\Db{\rm mod}$. Let $S$ be a simple module over $eRe$. Then  $S\otimesL_{eRe} eR\in \Db{R}$. Hence, $\Tor_m^{eRe}(S, eR)=0$ for sufficiently large $m$. Thus, $eR$ must have finite projective dimension. So $(a)$ holds.
	
	By \cite[Lemma 2.7]{AKLY2017}, the functor $\rHom_{eRe}(eR, -):{} \D{\Modcat{eRe}}\ra \D{\Modcat{R}}$ restricts to $\Db{\rm mod}$ if and only if the functor $\rHom_R(Re, -):{} \D{\Modcat{R}}\ra \D{\Modcat{eRe}}$ restricts to $\Kb{\rm proj}$. The functor $ \rHom_R(Re, -)$ restricts to $\Kb{\rm proj}$ if and only if $$eR\cong\Hom_R(Re, R)\cong\rHom_R(Re, R)\in \Kb{\pmodcat{eRe}}.$$
\end{proof}

\subsubsection{Existence of half recollements of bounded derived categories}

It is well known that the stratifying recollements (\ref{eq:7}) and (\ref{dual of 7}) restrict to recollements of bounded derived categories whenever $eRe$ has finite global dimension. Miyachi showed that this extends to the case where $eRe$ may have infinite global dimension, under the weaker assumption that $Re$ and $eR$ have finite projective dimension as $eRe$-modules. In this part, we revisit this result from a different perspective. Our approach differs from Miyachi's: rather than working directly with full recollements of bounded derived categories, we show that each of these conditions yields a half recollement of bounded derived categories.

\begin{Theo}\label{standardly}Let $R$ be a finite dimensional $k$-algebra, and $ReR$ be a stratifying ideal of $R$. Then: 
	\begin{enumerate}[label=(\alph*)]
		\item
		The condition $\pd_{eRe}(eR)<\infty$ holds if and only if one of the following equivalent statements holds:
		\begin{enumerate}[label=(a\arabic*)]
			\item
			There is a lower half stratifying recollement
			\begin{equation} \label{eq:hd}
				\xymatrix@=3.5em{
					\Db{{R/ReR}}\ar[rr]^{i_*}&&\Db{{R}}\ar[rr]^{\rHom_{R}(Re, -)}\ar@/^2pc/[ll]^{}&&\Db{{eRe}}\ar@/^2pc/[ll]^{}.}\end{equation}
			\item
			There is an upper half stratifying recollement
			\begin{equation}\label{eq:hu}
				\xymatrix@=3.2em{
					\Db{(R/ReR)^{\opp}}\ar[rr]^{i_*}&&\Db{R^{\opp}}\ar@/_2pc/[ll]^{}\ar[rr]^{\rHom_{R^{\opp}}(eR, -)}&&\Db{(eRe)^{\opp}}\ar@/_2pc/[ll]_{-\otimesL_{eRe}eR}.}\end{equation}\end{enumerate}

		\item
		The condition $\pd_{(eRe)^{\opp}}(Re)<\infty$ holds if and only if one of the following equivalent assertions holds:
		\begin{enumerate}[label=(b\arabic*)]
			\item
			There is an upper half stratifying recollement
			\begin{equation}
				\xymatrix@=3.5em{
					\Db{R/ReR}\ar[rr]^{i_*}&&\Db{R}\ar@/_2pc/[ll]^{}\ar[rr]^{\rHom_R(Re, -)}&&\Db{eRe}\ar@/_2pc/[ll]_{Re\otimesL_{eRe}-}.}\end{equation}
			\item
			There is a lower half stratifying recollement
			\begin{equation} 
				\xymatrix@=3.2em{
					\Db{(R/ReR)^{\opp}}\ar[rr]^{i_*}&&\Db{{R}^{\opp}}\ar[rr]^{\rHom_{R^{\opp}}(eR, -)}\ar@/^2pc/[ll]^{}&&\Db{(eRe)^{\opp}}\ar@/^2pc/[ll]^{}.}\end{equation}\end{enumerate}
	\end{enumerate}
\end{Theo}

\begin{proof}
	Assume that the half recollement (\ref{eq:hd}) or the half recollement (\ref{eq:hu}) exists. Then by Lemma \ref{lemmanew2}, $\pd_{eRe}(eR)<\infty$ holds. Conversely, assume that $\pd_{eRe}(eR)<\infty$ holds.
	
	Note that $$\Tor_m^{R/ReR}(R/ReR, M)=0=\Tor_m^{R/ReR}(N, R/ReR)$$ for all $M\in \modcat{R/ReR}$ and $N\in \modcat{(R/ReR)^{\opp}}$, $m\geq 1$. Then both $R/ReR\otimesL_{R/ReR}-$ and $-\otimesL_{R/ReR} R/ReR$ restrict to $\Db{\rm mod}$. So the functors $i_*$ always restrict to $\Db{\rm mod}$. By Lemma \ref{lemmanew2}, the functors $-\otimesL_{eRe} eR$ and $\rHom_{eRe}(eR, -)$ restrict to $\Db{\rm mod}$. Since $Re$ and $eR$ are projective as $R$-modules, the functors $Re\otimesL_{eRe} -$ and $-\otimesL_{eRe} eR$ restrict to $\Kb{\rm proj}$ while the functors $\rHom_R(Re, -)$ and $\rHom_{R^{\opp}}(eR, -)$ restrict to $\Db{\rm mod}$.
	
	Since $ReR$ is a stratifying ideal and $\pd_{eRe} eR<\infty$, $ReR\cong Re\otimesL_{eRe}eR\cong Re \otimes_{eRe} p(eR)\in \Kb{\pmodcat{R}}$, where $p(eR)$ is a finite projective resolution of $eR$ as a left $eRe$-module. Thus, $R/ReR\in \Kb{\pmodcat{R}}$ and therefore $-\otimesL_R R/ReR$ restricts to $\Db{\rm mod}$. Also since $R/ReR$ has finite projective dimension as left $R$-module, the functor $i_*=R/ReR\otimesL_{R/ReR}-$ restricts to $\Kb{\rm proj}$. By \cite[Lemma 2.7]{AKLY2017}, the right adjoint of $i_*$, $\rHom_R(R/ReR, -)$, restricts to $\Db{\rm mod}$. Analogously, the conditions (b), (b1) and (b2) are equivalent.
\end{proof}

As a corollary, we reobtain \cite[Theorem 2.7]{ParshallScott1987} as well as the particular case of \cite[Theorem 1]{Miyachi1992} for stratifying ideals.

\begin{Koro}\label{bounded}Let $R$ be a finite dimensional $k$-algebra, and $ReR$ be a stratifying ideal of $R$. Then the following statements are equivalent. 
\begin{enumerate}[label=(\alph*)]
	\item
	$\pd_{eRe}(eR)<\infty$ and $\pd_{(eRe)^{\opp}}(Re)<\infty$.
	\item
	There is a stratifying recollement
	\begin{equation}
		\xymatrix@=3.5em{
			\Db{R/ReR}\ar[rr]^{i_*}&&\Db{R}\ar@/^2pc/[ll]^{}\ar@/_2pc/[ll]_{}\ar[rr]^{\rHom_R(Re, -)}&&\Db{eRe}\ar@/^2pc/[ll]_{}\ar@/_2pc/[ll]_{Re\otimesL_{eRe}-}.}\end{equation}
	\item
	There is a stratifying recollement
	\begin{equation}
		\xymatrix@=3.2em{
			\Db{(R/ReR)^{\opp}}\ar[rr]^{i_*}&&\Db{R^{\opp}}\ar@/^2pc/[ll]^{}\ar@/_2pc/[ll]_{}\ar[rr]^{\rHom_R(eR, -)}&&\Db{(eRe)^{\opp}}\ar@/^2pc/[ll]_{}\ar@/_2pc/[ll]_{-\otimesL_{eRe}eR}.}\end{equation}
\end{enumerate}
\end{Koro}

Stratifying ideals satisfying the condition (a) of Corollary \ref{bounded} are recollement ideals in the sense of \cite{Miyachi1992}.
In particular, we re-obtain \cite[Lemma 3.4]{AKL2011} as a particular case.

\begin{Koro} Let $R$ be a finite dimensional algebra over $k$. Let $ReR$ be a stratifying ideal of $R$ such that $eRe$ has finite global dimension. Then there are (stratifying) recollements
\begin{align*}
	\xymatrix@=3.5em{
		\Db{R/ReR}\ar[rr]^{i_*}&&\Db{R}\ar@/^2pc/[ll]^{}\ar@/_2pc/[ll]_{}\ar[rr]^{\rHom_R(Re, -)}&&\Db{eRe}\ar@/^2pc/[ll]_{}\ar@/_2pc/[ll]_{Re\otimesL_{eRe}-},}\\\xymatrix@=3.2em{
		\Db{(R/ReR)^{\opp}}\ar[rr]^{i_*}&&\Db{R^{\opp}}\ar@/^2pc/[ll]^{}\ar@/_2pc/[ll]_{}\ar[rr]^{\rHom_{R^{\opp}}(eR, -)}&&\Db{(eRe)^{\opp}}\ar@/^2pc/[ll]_{}\ar@/_2pc/[ll]_{-\otimesL_{eRe}eR}.}
\end{align*}
\end{Koro}

\subsection{Characterisations of standardly stratified algebras via recollements} \label{subsection-characterisations}

In this subsection, we combine the results on stratifying recollements and standardly stratified algebras to prove the main result. As an immediate corollary, we derive the properly stratified analogue. More precisely, the bounded derived category of a properly stratified algebra admits an inductive characterisation from the bounded derived categories of local algebras.

\begin{proof}[Proof of Theorem \ref{equi}] Fix $e=e_n$.
  ``$(a2)\Rightarrow (a1)$" Assume that (a2) holds. By
  Remark \ref{progen} and Corollary \ref{coro1}, $ReR$ is a stratifying ideal
  and $ReR\cong Re\otimes_{eRe}eR\in \add_R Re$ since $eR\in \add_{eRe}eRe$. 
	Since $\Delta(n)=Re_n$, we have $Re_nR\in \mathcal{F}(\Delta(n))$. By (a2), $R/Re_nR\in \mathcal{F}(\Delta(1), \ldots, \Delta(n-1))$, and $\mathcal{F}(\Delta)$ is closed under extensions. Then the exact sequence
	$$0\lra Re_nR\lra R\lra R/Re_nR\lra 0$$
	in $\modcat{R}$ implies that
  ${}_RR\in\mathcal{F}(\Delta)$. Therefore, $R$ is left standardly stratified.   
	\medskip

        ``$(a3)\Leftrightarrow(a4)$" 
        If one of the recollements (\ref{d}) or
        (\ref{u}) exists, then at least one of the functors
        ${R/ReR \otimesL_{R/ReR} - }$
        or ${- \otimesL_{R/ReR}R/ReR}$ is fully faithful. 
        By Definition \ref{stratifyingidealdef}, $ReR$ is stratifying. By Lemma
        \ref{lemmanew2}, $\pd_{eRe}(eR) < \infty$. Thus
        Theorem \ref{standardly} implies that (a3) and (a4) are equivalent.

        \medskip

        ``$(a4)\Rightarrow (a2)$" As just shown in ``$(a3)\Leftrightarrow(a4)$",
        $Re_nR$ is stratifying and $\pd_{eRe}(eR) < \infty$. The idempotent
        $e$ is primitive. Thus $eR$ is projective as a left $eRe$-module 
        since the finitistic dimension $\findim (eRe)$ vanishes.
        By Corollary \ref{coro1}, $ReR$ is projective as
        a left $R$-module. Hence, (a2) holds.

        \medskip
        
        ``$(a1)\Rightarrow (a4)$" This implication follows from Lemma
        \ref{CM2021} and Theorem \ref{standardly}.
\end{proof}

\begin{Koro}\label{equiright}
	Let $(R, \{e_1, e_2, \ldots, e_n\})$ be a finite dimensional $k$-algebra, and $\{e_i\mid 1\leqslant i\leqslant n\}$ be a complete set of primitive pairwise orthogonal idempotents. Denote by $\overline{e}_i$ the image of $e_i$ in $R/Re_nR$. Then:
	\begin{enumerate}
		\item[(b)] The following are equivalent:
		\begin{enumerate}[label=(b\arabic*)]
			\item
			$(R, \{e_1, e_2, \ldots, e_n\})$ is a properly stratified algebra with $1 < \ldots < n$.
			\item
			$(R/Re_nR, \{\overline{e}_1, \overline{e}_2, \ldots, \overline{e}_{n-1}\})$ is a properly stratified algebra, and $Re_nR$ is projective both as a left and as a right $R$-module.
			\item
			$(R/Re_nR, \{\overline{e}_1, \overline{e}_2, \ldots, \overline{e}_{n-1}\})$ is a properly stratified algebra,
			and there is a recollement
			\begin{equation}
				\xymatrix@C=3.5em{
					\Db{{R/Re_nR}}\ar[rr]^{R/Re_nR\otimesL_{R/Re_nR}-}&&\Db{{R}}\ar@/^2pc/[ll]^{}\ar@/_2pc/[ll]_{}\ar[rr]^{\rHom_R(Re_n, -)}&&\Db{{e_nRe_n}}\ar@/^2pc/[ll]_{}\ar@/_2pc/[ll]_{Re_n\otimesL_{e_nRe_n}-}.}\label{eq22}\end{equation}
			\item
			$(R/Re_nR, \{\overline{e}_1, \overline{e}_2, \ldots, \overline{e}_{n-1}\})$ is a properly stratified algebra,
			and there is a recollement
			\begin{equation}
				\xymatrix@C=3.5em{
					\Db{(R/Re_nR)^{\opp}}\ar[rr]^{-\otimesL_{Re_nR}R/Re_nR}&&\Db{{R^{\opp}}}\ar@/^2pc/[ll]^{}\ar@/_2pc/[ll]_{}\ar[rr]^{\rHom_{R^{\opp}}(e_nR, -)}&&\Db{(e_nRe_n)^{\opp}}\ar@/^2pc/[ll]_{}\ar@/_2pc/[ll]_{-\otimesL_{e_nRe_n}e_nR}.}\end{equation}
			
		\end{enumerate}
		
		\item[(c)]
		The following are equivalent: 
		
		\begin{enumerate}[label=(c\arabic*)]
			\item
			$(R, \{e_1, \ldots, e_n\})$ is a right standardly stratified algebra with $1 < \cdots < n$.
			\item
			$(R/Re_nR, \{\overline{e}_1, \ldots, \overline{e}_{n-1}\})$ is a right standardly stratified algebra, and $Re_nR$ is a projective right $R$-module.
			\item
			$(R/Re_nR, \{\overline{e}_1, \ldots, \overline{e}_{n-1}\})$ is a right standardly stratified algebra, and there is an upper half stratifying recollement
			\begin{equation} 
				\xymatrix@C=3.5em{
					\Db{{R/Re_nR}}\ar[rr]_{R/Re_nR\otimesL_{R/Re_nR}-}&&\Db{{R}}\ar[rr]_{\rHom_R(Re_n, -)}\ar@/_2pc/[ll]_{R/Re_nR\otimesL_R-}&&\Db{{e_nRe_n}}\ar@/_2pc/[ll]_{Re_n\otimesL_{e_nRe_n}-}.}\end{equation}
			\item
			$(R/Re_nR, \{\overline{e}_1, \ldots, \overline{e}_{n-1}\})$ is a right standardly stratified algebra, and there is a lower half stratifying recollement
			
			\begin{equation} 
				\xymatrix@C=3.5em{
					\Db{{R/Re_nR}}\ar[rr]^{R/Re_nR\otimesL_{R/Re_nR}-}&&\Db{{R}}\ar[rr]^{\rHom_R(Re_n, -)}\ar@/^2pc/[ll]^{\rHom_R(R/Re_nR, -)}&&\Db{{e_nRe_n}}\ar@/^2pc/[ll]^{\rHom_{e_nRe_n}(e_nR, -)}.}\end{equation}
			
		\end{enumerate} 
	\end{enumerate} \label{equi2}
\end{Koro}
\begin{proof}
(c) follows since an algebra is right standardly stratified if and only if its opposite algebra is left standardly stratified with respect to the same ordering.
		
		By definition, the statement (b$i$) holds if and only if the statements (a$i$) of Theorem \ref{equi} and (c$i$) hold for every $i=1, \ldots, 4$. Hence (b) follows from Theorem \ref{equi}  and (c).
\end{proof}

Thanks to the main theorem, we can find numerous non-trivial examples of half recollements of bounded derived categories. The following algebra is a classic example (see for instance \cite[Page 42]{Cline1996}) of a standardly stratified algebra which is not properly stratified. Moreover, the following example presents an algebra in which every idempotent ideal arising from a primitive idempotent is stratifying (see also \cite[Example 3]{Auslander1992b}).
In particular, the bounded derived category of the algebra below fits into non-trivial lower half and upper half recollements of bounded derived categories.

\begin{Example}
	Let $R$ be the quotient of the path algebra of the quiver  $$\xymatrix{
		1\ar@/^1pc/[rr]^{\alpha}&&2\ar@/^1pc/[ll]^{\beta}
	}$$ modulo the relation $\beta\alpha\beta=0$. Denote by $e_i$ the primitive idempotent associated with the vertex $i$ for $i=1, 2$. $R$ can be written as 
	$$\begin{pmatrix}
		{\rm span}\{e_1, \beta \alpha \}&{\rm span}\{\beta\}\\
		{\rm span}\{\alpha, \alpha\beta\alpha\}&{\rm span}\{e_2, \alpha\beta\}
	\end{pmatrix}.
	$$
	Hence the dimension of $R$ is $7$, where $Re_1$ has dimension $4$ and $Re_2$ has dimension 3. $R$ is left standardly stratifying with $\Delta(1)=1$ and $\Delta(2)=Re_2$. Hence, $Re_2R$ is a stratifying ideal and projective as left $R$-module. As $\overline{\Delta}(1)=1$ and $\overline{\Delta}(2)=\indt{2}{1}{}$, $R$ is not properly stratified, and in particular it is not right standardly stratified with respect to this ordering. By Theorem \ref{equi} and Corollary \ref{equi2}, there are half recollements
	\begin{align*}
		\xymatrix@=3.5em{
			\Db{{R/Re_2R}}\ar[rr]^{R/Re_2R\otimesL_{R/Re_2R}-}&&\Db{{R}}\ar[rr]^{\rHom_R(Re_2, -)}\ar@/^2pc/[ll]^{}&&\Db{{e_2Re_2}}\ar@/^2pc/[ll]^{},} \\ \xymatrix@=3.5em{
			\Db{(R/Re_2R)^{\opp}}\ar[rr]^{-\otimesL_{Re_2R}R/Re_2R}&&\Db{{R^{\opp}}}\ar@/_2pc/[ll]_{}\ar[rr]^{\rHom_{R^{\opp}}(e_2R, -)}&&\Db{(e_2Re_2)^{\opp}}\ar@/_2pc/[ll]_{-\otimesL_{e_2Re_2}e_2R}.} 
	\end{align*}
	
	Now we consider the idempotent ideal $Re_1R$ whose basis is $\{e_1, \alpha, \beta, \alpha\beta, \beta\alpha, \alpha\beta\alpha\}$. We claim that $Re_1R$ is also a stratifying ideal. We shall show it directly without using Theorem \ref{equi}. 
	The basis of $e_2Re_1$ is $\{\alpha, \alpha\beta\alpha\}$, and $e_2Re_1$ is isomorphic to $e_1Re_1$ as right $e_1Re_1$-modules. So $Re_1$ is a projective right $e_1Re_1$-module. It is straightforward to check the multiplication
        ${\mu:{}Re_1\otimes_{e_1Re_1}e_1R\ra Re_1R}$ is an isomorphism of $R$-$R$-bimodules. Hence, $Re_1R$ is a stratifying ideal of $R$. 
	Since $Re_1R$ has dimension $6$, it is not projective as left $R$-module. By Corollary \ref{coro1}, it follows that $e_1R$ has infinite projective dimension as $e_1Re_1$-module since $e_1Re_1$ is local and $\findim (e_1Re_1)=0$. Indeed, $e_1Re_1=k[x]/(x^2)$ and $e_1R$ has dimension $3$.
	Since $Re_1R$ is a stratifying ideal of $R$, there are recollements
	\begin{equation}\label{l}
		\xymatrix@C=3.5em{
			\D{\Modcat{R/Re_1R}}\ar[rr]^{R/Re_1R\otimesL_{R/Re_1R}-}&&\D{\Modcat{R}}\ar@/^2pc/[ll]^{}\ar@/_2pc/[ll]_{}\ar[rr]^{\rHom_R(Re_1, -)}&&\D{\Modcat{e_1Re_1}}\ar@/^2pc/[ll]_{}\ar@/_2pc/[ll]_{Re_1\otimesL_{e_1Re_1}-},}\end{equation}
	\begin{equation}\label{r}
		\xymatrix@C=3.5em{
			\D{\rModcat{R/Re_1R}}\ar[rr]^{-\otimesL_{R/Re_1R}R/Re_1R}&&\D{\rModcat{R}}\ar@/^2pc/[ll]^{}\ar@/_2pc/[ll]_{}\ar[rr]^{\rHom_{R^{\opp}}(e_1R, -)}&&\D{\rModcat{e_1Re_1}}\ar@/^2pc/[ll]_{}\ar@/_2pc/[ll]_{-\otimesL_{e_1Re_1}e_1R}.}\end{equation}
	By Corollary \ref{bounded}, neither (\ref{l}) nor (\ref{r}) restricts to $\Db{{\rm mod}}$. Actually,  we get from (\ref{l}) and (\ref{r}) the following two half recollements 
	\begin{align*}
		\xymatrix@C=3.5em{
			\Db{{R/Re_1R}}\ar[rr]^{R/Re_1R\otimesL_{R/Re_1R}-}&&\Db{{R}}\ar@/_2pc/[ll]_{}\ar[rr]^{\rHom_R(Re_1, -)}&&\Db{{e_1Re_1}}\ar@/_2pc/[ll]_{Re_1\otimesL_{e_1Re_1}-},} \\
		\xymatrix@=3.2em{
			\Db{{(R/Re_1R)^{\opp}}}\ar[rr]^{-\otimesL_{R/Re_1R}R/Re_1R}&&\Db{{R^{\opp}}}\ar@/^2pc/[ll]_{}\ar[rr]^{\rHom_{R^{\opp}}(e_1R, -)}&&\Db{(e_1Re_1)^{\opp}}\ar@/^2pc/[ll]}.
	\end{align*}
	Since $R/Re_1R$ is isomorphic to $k$, it follows by Corollary \ref{equi2} that $(R, \{f_1, f_2\})$ is a right standardly stratified algebra with $f_1=e_2$ and $f_2=e_1$, but it is not properly stratified.
\end{Example}

\medskip

\section{Gorenstein standardly stratified algebras}

A finite dimensional $k$-algebra $R$ is called {\em Iwanaga--Gorenstein} if $\id_R(R)<\infty$ and $\id_{R^{\opp}}(R)<\infty$. 

The main aim of this section is to prove the following result:

\begin{Theo}\label{mainresult}
  Let $R$ be an Iwanaga--Gorenstein $k$-algebra. Then the
  following assertions are equivalent.
	\begin{enumerate}[label=(\alph*)]
		\item $(R, \{e_1, \ldots, e_n\})$ is  a properly stratified algebra with the ordering $1 < 2 < \cdots < n$.
		\item  $(R, \{e_1, \ldots, e_n\})$ is  a left standardly stratified algebra with the ordering $1 < 2 < \cdots < n$.
		\item $(R, \{e_1, \ldots, e_n\})$ is  a right standardly stratified algebra with the ordering $1 < 2 < \cdots < n$.
	\end{enumerate}
\end{Theo}

Recollements serve as the foundational structure to understand and prove
this fact. Indeed, the proof is built on two key ideas:

\begin{itemize}
	\item If a given left standardly stratified algebra $(R, \{e_1, \ldots, e_n\})$ is Iwanaga--Gorenstein, then the algebras $e_nRe_n$ and $R/Re_nR$ are also Iwanaga--Gorenstein. (This fact is a consequence of Lemma \ref{lssaproperties}).
	\item As evidenced by \cite{CK2017}, recollements of bounded derived categories of finite dimensional algebras can detect Gorensteinness of the algebras involved and when all algebras are Iwanaga--Gorenstein, the recollement restricts to one of bounded homotopy categories of projectives (resp. of injectives).
	As it turns out the same occurs for half recollements of bounded derived categories of finite dimensional algebras. In particular, every half recollement of bounded derived categories of Iwanaga--Gorenstein algebras is a full recollement (see Proposition \ref{halfrecollementsGorenstein} and Proposition \ref{GR}).
\end{itemize}

Once these two key facts are established, Theorem \ref{mainresult} will follow from Theorem \ref{equi}. Indeed, a left standardly stratified algebras gives rise to a half recollement of bounded derived categories. By the facts above, such a half recollement necessarily extends to a full recollement of bounded derived categories. Consequently, Corollary \ref{equi2} implies that the original algebra must be properly stratified. Details of the proof will be given after
Lemma \ref{lssaproperties} has been stated and shown.

\subsection{Gorensteinness and recollements}

Recollements can be viewed as analogues of short exact sequences between triangulated categories and their existence can be used to transfer homological properties between different derived categories of finite dimensional algebras, like finiteness of global dimension or more generally Gorensteinness.
The following connection between recollements and Gorensteinness extends \cite[Proposition 3.4]{CK2017}.

\begin{Prop}\label{GR} 
Let $R, R_1, R_2$ be finite dimensional $k$-algebras. 
Suppose there is a recollement
$$\xymatrix{
  \Db{R_1}\ar[rr]^{i_*=i_!}&&\Db{R}\ar@/^1.5pc/[ll]^{i^!}\ar@/_1.5pc/[ll]_{i^*}\ar[rr]^{j^!=j^*}&&\Db{R_2}\ar@/_1.5pc/[ll]_{j_!}\ar@/^1.5pc/[ll]^{j_*}}
\phantom{xxxxx} (\ast).$$
Then the following statements are equivalent:
\begin{enumerate}[label=(\alph*)]
\item
The algebra $R$ is Iwanaga--Gorenstein.
\item
The algebras $R$, $R_1$ and $R_2$ are Iwanaga--Gorenstein.
\item
The algebras $R_1$ and $R_2$ are Iwanaga--Gorenstein and moreover
$i^!(R)\in\Kb{\pmodcat{R_1}}$.
\item
The algebras $R_1$ and $R_2$ are Iwanaga--Gorenstein and moreover the given 
recollement  $(\ast)$
restricts to a recollement of bounded homotopy categories of 
projectives.
\item The algebras $R_1$ and $R_2$ are Iwanaga--Gorenstein and moreover the
  given recollement $(\ast)$
  restricts to a recollement of bounded homotopy categories of 
injectives.
\item The algebras $R_1$ and $R_2$ are Iwanaga--Gorenstein and $j_!$ restricts to $\Kb{\imodcat{R_2}}\rightarrow \Kb{\imodcat{R}}$.
\item The algebras $R_1$ and $R_2$ are Iwanaga--Gorenstein and $j_!j^!(R)\in \Kb{\imodcat{R}}$.
\end{enumerate}
\end{Prop}
\begin{proof}
	The equivalences (a) $\Leftrightarrow$ (b) $\Leftrightarrow$ (c) $\Leftrightarrow$ (d) are \cite[Proposition 3.4]{CK2017}. 
	
	Assume that (a) holds. Since the assertion (a) is equivalent to (d) and $\Kb{\operatorname{proj}}=\Kb{\operatorname{inj}}$ for all the algebras involved, the assertion (e) follows. The implication (e) $\Rightarrow$ (f) is immediate.
	
Assume that (f) holds. Since $j_*$ restricts to $\Db{{\rm mod}}$ the functor $j^*$ restricts to $\Kb{\operatorname{proj}}$ by Lemma \ref{adjoint}. Thus $j^*(R)\in \Kb{\pmodcat{R_2}}=\Kb{\imodcat{R_2}}$ since $R_2$ is Iwanaga--Gorenstein. By (f), $j_{!}j^*(R)\in \Kb{\imodcat{R}}$. So (g) holds.

It remains to show the implication (g) $\Rightarrow$ (a). To do that, we show that $\id_R(R)<\infty$ and $\pd_R (D(R))<\infty$. Assume that (g) holds.
Consider the canonical triangle \begin{align}
	j_!j^!(R)\rightarrow R\rightarrow i_*i^*(R)\rightarrow j_!j^!(R)[1]. \label{eq29}
\end{align} By (g), $j_!j^!(R)\in \Kb{\imodcat{R}}$. 
By Lemma \ref{adjoint}, the functors $i^*, i_*, j_!, j^*$ restrict to $\Kb{\operatorname{proj}}$ while by Lemma \ref{adjointdual}, the functors $i^!, i_*, j^*, j_*$ restrict to $\Kb{\operatorname{inj}}$.
So, $i^*(R)\in \Kb{\pmodcat{R_1}}=\Kb{\imodcat{R_1}}$ since $R_1$ is Iwanaga--Gorenstein. Thus, $i_*i^*(R)\in \Kb{\imodcat{R}}$ and (\ref{eq29}) implies that ${R\in \Kb{\imodcat{R}}}$. This shows that $\id_R(R)<\infty$. 

To show that $\pd_R (D(R))<\infty$, consider the canonical triangle
\begin{align}
	j_!j^!(D(R))\rightarrow D(R)\rightarrow i_*i^*(D(R))\rightarrow j_!j^!(D(R))[1]. \label{eq30}
\end{align}
We claim that $j_!$ and $i^*$ restrict to $\Kb{\operatorname{inj}}$. Since $R_2$ is Iwanaga--Gorenstein, $D(R_2)\in \Kb{\pmodcat{R_2}}$ and $j_!(D(R_2))\in \Kb{\pmodcat{R}}$. By (g), we obtain that $j_!(D(R_2))\cong j_!j^!(j_!(D(R_2)))\in \Kb{\imodcat{R}}$.
It follows that $j_!j^!$ restricts to $\Kb{\operatorname{inj}}$ and thus $j_!j^!(D(R))\in \Kb{\imodcat{R}}$. By (\ref{eq30}), we deduce that $ i_*i^*(D(R))\in \Kb{\imodcat{R}}$.
Thus, $i^*(D(R))\cong i^!i_*i^*(D(R))\in \Kb{\imodcat{R_1}}$ because $i^!$ restricts to $\Kb{\operatorname{inj}}$. So, the claim follows. 

As $R_1$ is Iwanaga--Gorenstein, $i^*(D(R))\in \Kb{\pmodcat{R_1}}$, and thus, $i_*i^*(D(R))\in \Kb{\pmodcat{R}}$. Since $R_2$ is Iwanaga--Gorenstein, $j^!(D(R))\in \Kb{\pmodcat{R_2}}$ and thus $j_!j^!(D(R))\in \Kb{\pmodcat{R}}$. By (\ref{eq30}), $D(R)\in \Kb{\pmodcat{R}}$. So, (a) holds.
\end{proof}

Can we replace the recollement in Proposition \ref{GR}  by a half recollement?
Perhaps surprisingly, the answer is essentially affirmative: half recollements between bounded derived categories of Iwanaga--Gorenstein algebras which descend from $\D{\operatorname{Mod}}$-recollements can, in fact, be extended to full recollements of bounded derived categories.

\begin{Prop}\label{halfrecollementsGorenstein}Let $R, R_1, R_2$ be Iwanaga--Gorenstein algebras. Suppose there is a recollement 
$$\xymatrix{
	\D{\Modcat{R_1}}\ar[rr]^{i_*=i_!}&&\D{\Modcat{R}}\ar@/^1.5pc/[ll]^{i^!}\ar@/_1.5pc/[ll]_{i^*}\ar[rr]^{j^!=j^*}&&\D{\Modcat{R_2}}\ar@/_1.5pc/[ll]_{j_!}\ar@/^1.5pc/[ll]^{j_*}}.$$
such that the functors $i_!$, $i^!$, $j^*$ and $j_*$ restrict to $\Db{\operatorname{mod}}$.
Then there is a recollement $$\xymatrix{
	\Db{R_1}\ar[rr]^{i_*=i_!}&&\Db{R}\ar@/^1.5pc/[ll]^{i^!}\ar@/_1.5pc/[ll]_{i^*}\ar[rr]^{j^!=j^*}&&\Db{R_2}\ar@/_1.5pc/[ll]_{j_!}\ar@/^1.5pc/[ll]^{j_*}}.$$
\end{Prop}
\begin{proof}
	By Lemma \ref{adjoint}, the functors $i_!$ and $j^*$ restrict to $\Kb{\operatorname{proj}}$. Since $R, R_1$ and $R_2$ are Iwanaga--Gorenstein, $\Kb{\pmodcat{R}}=\Kb{\imodcat{R}}$ and $\Kb{\pmodcat{R_i}}=\Kb{\imodcat{R_i}}$ for $i=1, 2$. Thus, the functors $i_!$ and $j^*$ restrict to $\Kb{\operatorname{inj}}$. By Lemma \ref{adjointdual}, the functors $i^*$ and $j_!$ restrict to $\Db{\operatorname{mod}}$ and the result follows.
\end{proof}

As an application of Proposition \ref{GR} to properly stratified algebras, the following corollary is now immediate. 

\begin{Koro}\label{GRP} Let $(R, \{e_1, \ldots, e_n\})$ be  a properly stratified algebra with ordering $1 < 2 < \ldots < n$. Then the following statements are equivalent. 
	\begin{enumerate}[label=(\alph*)]
		\item
		The algebra $R$ is Iwanaga--Gorenstein.
		\item
		The algebras $e_nRe_n$ and $R/Re_nR$ are Iwanaga--Gorenstein and $\rHom_R(R/Re_nR, R)$ is in $\Kb{\pmodcat{R/Re_nR}}$.
		\item
		The algebras $e_nRe_n$ and $R/Re_nR$ are Iwanaga--Gorenstein and moreover the 
		recollement (\ref{eq22}) restricts to a recollement of bounded homotopy categories of 
		projectives.
		\item The algebras $e_nRe_n$ and $R/Re_nR$ are Iwanaga--Gorenstein and moreover the  
		recollement (\ref{eq22}) restricts to a recollement of bounded homotopy categories of 
		injectives.
		\item The algebra $R/Re_nR$ is Iwanaga--Gorenstein and $\id_R(Re_n)<\infty$. \\ In particular, $j_!j^!(R)=Re_nR\in \Kb{\imodcat{R}}$.
	\end{enumerate}
\end{Koro}
\begin{proof} By Theorem \ref{equi} and Corollary \ref{equi2}, there exists a recollement between the bounded derived categories of $R/Re_nR$, $R$ and $e_nRe_n$. Thus, the statements (a), (b), (c), (d) are equivalent by Proposition \ref{GR}.
	Write $e=e_n$. 
	Observe that $j_!j^!(R)=Re\otimesL_{eRe} eR\otimesL_R R\cong Re\otimesL_{eRe} eR\cong ReR$. Moreover, $\add ReR=\add Re$ by  Lemma \ref{CM2021}. Thus, by Proposition \ref{GR}(g), the statement (e) holds.

        Assume now that (e) holds. 
	By assumption, $\pd_{R^{\opp}}D(Re)<\infty$. So there exists a projective resolution $0\rightarrow P_t\rightarrow \cdots\rightarrow P_0\rightarrow D(Re)\rightarrow 0$. Applying the exact functor $-\otimes_R Re$ we get an exact sequence $0\rightarrow P_te\rightarrow \cdots \rightarrow P_0e\rightarrow D(eRe)\rightarrow 0$ and each $P_ie\in \add Re$. By Corollary \ref{equi2}, $Re$ is projective as right $eRe$-module and thus $D(eRe)$ has finite projective dimension as right $eRe$-module. Since $eRe$ is local, its finitistic dimension is zero and thus $D(eRe)$ is projective as right $eRe$-module. Therefore, $eRe$ is self-injective. So assertion (g) in Proposition \ref{GR} is satisfied and thus $R$ is Iwanaga--Gorenstein.  So (a) holds.
\end{proof}

This corollary also illustrates why extending Proposition 3.4 of \cite{CK2017} to Proposition 3.1 is worthwhile. In the context of properly stratified algebras, the object $j_!j^!(R)$ can be described explicitly. Indeed, it is the module $ReR$ with $e=e_n$. In contrast, it is less clear what is the object $i^!(R)$ in this context. Corollary \ref{GRP} thus complements the main result of \cite{CM2021} by providing additional tools to test whether a properly stratified algebra is Iwanaga--Gorenstein.

We now turn our attention to left (resp. right) standardly stratified algebras. Before proving the main result, we give some properties of Iwanaga--Gorenstein standardly stratified algebras.

\begin{Lem}\label{lssaproperties}
	 Let $(R, \{e_1, \ldots, e_n\})$ be  left standardly stratified with the ordering $1 < 2 < \cdots < n$.
	 \begin{enumerate}[label=(\alph*)]
	 	\item If $\id_R(R)<\infty$, then $R$ is Iwanaga--Gorenstein.
	 	\item If $R$ is  Iwanaga--Gorenstein, then $e_n R e_n$ is self-injective.
	 	\item The following assertions are equivalent:
	 	\begin{enumerate}[label=(c\arabic*)]
	 		\item $R$ is Iwanaga--Gorenstein.
	 		\item Every module in $\mathcal{F}(\Delta)$ has finite injective dimension.
	 		\item $\id_R(\Delta(i))<\infty$ for every $i=1, \ldots, n$.
	 	\end{enumerate}
	 \end{enumerate}
\end{Lem}
\begin{proof}
	By Theorem 3.1 and Theorem 2.1 of \cite{Agoston2000a}, $R$ and $R^{\opp}$ have finite finitistic dimension. Thus, by Proposition 1.2 of \cite{CM2021}, (a) follows.
	
	By Theorem 1.6, Proposition 2.2 and Lemma 2.5 of \cite{Agoston2000}, $\mathcal{F}(\Delta)$ is contravariantly finite, $\mathcal{F}(\Delta)\cap \mathcal{F}(\overline{\nabla})=\add T$, where $T$ is the characteristic tilting module. So the proof of Theorem 2.5(1) of \cite{CM2021} is also valid for left standardly stratified algebras, and so (b) holds.

        In any short exact sequence, if two out of three terms have finite injective dimension, then so does the third. It follows that (c2) and (c3) are equivalent.
	
	By definition, $R\in \mathcal{F}(\Delta)$. Thus, by (a), the implication (c2) $\Rightarrow$ (c1) is clear.
	
	Finally assume (c1) holds. We will prove by induction on $l$ that all modules in $\mathcal{F}(\{\Delta(1), \ldots, \Delta(l)\})$ have finite injective dimension. Since $R$ is Iwanaga--Gorenstein the characteristic tilting module $T$ is also cotilting (at the moment we do not claim that it is a characteristic cotilting module). So $T$ has finite injective dimension. The direct summand $\Delta(1)$ of $T$ has finite injective dimension and so all modules in $\mathcal{F}(\Delta(1))$ have finite injective dimension. 
	Assume that $l>1$. Consider the minimal $\mathcal{F}(\overline{\nabla})$-approximation $0\rightarrow \Delta(l)\rightarrow T(l)\rightarrow X(l)\rightarrow 0$. So $X(l)\in \mathcal{F}(\{ \Delta(1), \cdots, \Delta(l-1) \})$ and by induction $\id_R X(l)<\infty$. Since $R$ is Iwanaga--Gorenstein and $T(l)$ is a direct summand of the tilting-cotilting module $T$ we get that $\id_R(T(l))<\infty$. Thus $\id_R(\Delta(l))<\infty$. Therefore all modules $\mathcal{F}(\Delta(1), \ldots, \Delta(l))$ have finite injective dimension and (c2) holds.
\end{proof}

 We are now ready to prove the main result of this section.

\begin{proof}[Proof of Theorem \ref{mainresult}] By definition of properly stratified algebras, (a) implies both (b) and (c).
	Assume that (b) holds. 
	To prove (a) we will proceed by induction on $n$. Assume that $n=1$. In this case, $R$ is a local algebra and it is automatic that $R$ is a properly stratified algebra. Assume that $n>1$. As $R$ is Iwanaga-Gorenstein by assumption,  Lemma \ref{lssaproperties} implies that $e_nRe_n$ is a self-injective algebra and every module in $\mathcal{F}(\Delta)$ has finite injective dimension. By Theorem \ref{equi}, $(R/Re_nR, \{\overline{e}_1, \overline{e}_2, \cdots, \overline{e}_{n-1}\})$ is a standardly stratified algebra and there is a lower half stratifying recollement
	\begin{equation} \label{eq31}
		\xymatrix@C=3.5em{
			\Db{{R/Re_nR}}\ar[rr]^{R/Re_nR\otimesL_{R/Re_nR}-}&&\Db{{R}}\ar[rr]^{\rHom_R(Re_n, -)}\ar@/^2pc/[ll]^{}&&\Db{{e_nRe_n}}\ar@/^2pc/[ll]^{}.}\end{equation}
		Thus, $Re_nR$ is a stratifying ideal (see Definition \ref{stratifyingidealdef}(d)). Hence, $R/Re_nR\in \mathcal{F}(\Delta)$ and consequently $\id_R(R/Re_nR)<\infty$ by  Lemma \ref{lssaproperties}. Thus, $\Ext_{R/Re_nR}^i(M, R/Re_nR)\cong \Ext_R^i(M, R/Re_nR)=0$ for all $M\in \modcat{R/Re_nR}$ and $i>\id_R(R/Re_nR)$. This means that $\id_{R/Re_nR}(R/Re_nR)<\infty$. By Lemma \ref{lssaproperties}(a), $R/Re_nR$ is Iwanaga--Gorenstein.  It follows then by Proposition \ref{halfrecollementsGorenstein}, Lemma \ref{main} and (\ref{eq31}) that there exists a recollement \begin{equation}
			\xymatrix@C=3.5em{
				\Db{{R/Re_nR}}\ar[rr]^{R/Re_nR\otimesL_{R/Re_nR}-}&&\Db{{R}}\ar@/^2pc/[ll]^{}\ar@/_2pc/[ll]_{}\ar[rr]^{\rHom_R(Re_n, -)}&&\Db{{e_nRe_n}}\ar@/^2pc/[ll]_{}\ar@/_2pc/[ll]_{Re_n\otimesL_{e_nRe_n}-}.}\end{equation}
			By Corollary \ref{equi2}(b), $(R, \{e_1, \ldots, e_n\})$ is a properly stratified algebra. So (a) holds.
			
			If $(R, \{e_1, \ldots, e_n\})$ is right standardly stratified, then $(R^{\opp}, \{e_1, \ldots, e_n\})$ is left standardly stratified. By the implication (b) $\Rightarrow$ (a), $(R^{\opp}, \{e_1, \ldots, e_n\})$ is properly stratified and so (a) holds.
\end{proof}

We summarise a list of applications that follow by using our methodology for Iwanaga--Gorenstein left standardly stratified algebras. By our main result, it suffices to present these for properly stratified algebras.

\begin{Koro}\label{3}
	Let $(R, \{e_1, \ldots, e_n\})$ be an Iwanaga--Gorenstein properly stratified algebra with the ordering $1 < 2 < \cdots < n$. Then the following statements hold.
	\begin{enumerate}[label=(\alph*)]
		\item The algebras $R_i:=R/R(e_{i+1}+\cdots+ e_n)R$ are Iwanaga--Gorenstein properly stratified for all $i=1, \ldots, n-1$.
		\item The algebras $\End_R(\Delta(i))$ are self-injective and $\Delta(i)\in \Kb{\pmodcat{R}}=\Kb{\imodcat{R}}$ for every $i=1, \ldots, n$.
		\item $Re_nR$ is projective as a right and as a left $R$-module. In particular, it has finite injective dimension as a right and as a left $R$-module.
		\item $\add e_nR=\add D(Re_n)=\add e_nRe_n$ as left $e_nRe_n$-modules 
	\end{enumerate}
\end{Koro}
\begin{proof}
	By Corollary \ref{GRP} and Corollary \ref{equi2}, $R/Re_nR$ is Iwanaga--Gorenstein properly stratified. By induction, it follows that all algebras $R_i$ are Iwanaga--Gorenstein properly stratified. 
	
	The fact that $\Delta(i)\in  \Kb{\pmodcat{R}}=\Kb{\imodcat{R}}$ for every $i=1, \ldots, n$ follows by Lemma \ref{lssaproperties}. The fact that the algebras $\End_R(\Delta(i))$ are self-injective was proved in \cite[Theorem 2.5]{CM2021}; it is also a consequence of Corollary \ref{GRP}. Indeed, $\End_{R}(\Delta(i))\cong \End_{R_i}(\Delta(i))$ is Iwanaga--Gorenstein and local by Corollary \ref{GRP} since $R_i$ is Iwanaga--Gorenstein properly stratified by (a). So (b) holds since Iwanaga--Gorenstein local algebras are self-injective.

	The statement (c) follows from $\Kb{\pmodcat{R}}=\Kb{\imodcat{R}}$ and Theorem \ref{equi}. By Lemma \ref{CM2021} and by $R$ and $R^{\opp}$ being left standardly stratified algebras, $e_nR$ is a projective generator of $e_nRe_n$ and $Re_n$ is a projective generator of $e_nRe_n$ as a right module. Since $e_nRe_n$ is self-injective $\add D(Re_n)=\add D(e_nRe_n)=\add e_nRe_n$, assertion
(d) holds.
\end{proof}

\section{Examples}

For standardly stratified structures, the choice of ordering on the primitive idempotents is crucial. The following example shows that for non-Iwanaga--Gorenstein algebras with stratifications by self-injective local rings, some orderings give a properly stratified structure, while others yield only a standardly stratified structure.

\begin{Example}
	Let $R$ be the quotient of the path $k$-algebra of the quiver 
	\[
	\xymatrix{
		1 \ar[r]^{\alpha} & 2 \ar@(ur,dr)[]^{\beta}
	}
	\] modulo the relations $\beta^2=\beta\alpha=0$. The projective indecomposable $R$-modules are then $\indt{1}{2}{}$ and $\indt{2}{2}{}$. Denote by $e_i$ the primitive idempotent associated with the vertex $i$ for $i=1, 2$. With respect to this choice of primitive idempotents, the standard modules are $\Delta(1)=1$ and $\Delta(2)=Re_2$, while the proper standard modules are $\overline{\Delta}(1)=1$ and $\overline{\Delta}(2)=2$. So, $(R, \{e_1, e_2\})$ is right standardly stratified but not properly stratified. By Theorem \ref{mainresult}, $R$ is not Iwanaga--Gorenstein. Now fix ${f_1=e_2}$ and ${f_2=e_1}$. With respect to this ordering of the idempotents, the standard modules are the indecomposable projective modules and $\overline{\Delta}(1)=Re_1$ while $\overline{\Delta}(2)=2$. This means that $(R, \{f_1, f_2\})$ is properly stratified. By Theorem \ref{equi}, there exists a recollement
	\begin{equation}
		\xymatrix{
			\Db{{k[x]/(x^2)}}\ar[rr]&&\Db{{R}}\ar@/^1.5pc/[ll]^{}\ar@/_1.5pc/[ll]_{}\ar[rr]&&\Db{{k}}\ar@/^1.5pc/[ll]_{}\ar@/_1.5pc/[ll].}\end{equation}
	Since $R$ is not Iwanaga--Gorenstein, this recollement does not restrict to $\Kb{proj}$ even though the corner algebras are Iwanaga--Gorenstein.
\end{Example}

\medskip

Example 2.7 of \cite{CM2021} seems to cause a problem. Fortunately, this
example is not accurate. Indeed, there the authors computed
$\Hom_R(R/ReR, R)$ and then applied Proposition 3.4 of \cite{CK2017}.
However, Proposition 3.4 of \cite{CK2017} requires instead $\rHom_R(R/ReR, R)$.
Our next example corrects example 2.7 of \cite{CM2021}.

\begin{Example} Let $R$ be the opposite algebra of the algebra considered in Example 2.7 of \cite{CM2021} (we use this $R$ since we are using left $R$-modules instead of right $R$-modules). Hence, $R$ is the quotient of the path $k$-algebra of the quiver 
	$\xymatrix@C=3cm{
		1 
		\ar@<0.7ex>[r]^{\alpha_1}
		& 
		2 
		\ar@<0.7ex>[r]^{\alpha_2}
		\ar@<0.7ex>[l]^{\beta_1}
		& 
		3 
		\ar@<0.7ex>[l]^{\beta_2}
	}
	$ modulo the relations $\beta_1\alpha_1\alpha_2=\beta_2\beta_1\alpha_1=\beta_2\alpha_2=0$, $\alpha_2\beta_2=\beta_1\alpha_1\beta_1\alpha_1$. Let $e_i$ be the primitive idempotent associated with the vertex $i$ for $i=1, 2, 3$.

	The  indecomposable projective $R$-modules of this algebra can be
        displayed in terms of their composition factors as follows:
	$$ Re_1=\begin{smallmatrix}
	&	1 & \\ &2 & \\ 1 & & 3 \\ 2 & & \\ 1 & & \\ & 2 & \\ & 1 &
	\end{smallmatrix}, \qquad Re_2=\begin{smallmatrix}
	& 2 & \\ 1 & & 3 \\ 2 & & \\ 1 & & \\ & 2 & \\ & 1 &  
\end{smallmatrix}, \qquad Re_3=\indt{3}{2}{1}. $$ We shall compute $\rHom_R(R/Re_3R, R)$. 
Let $I(i)$ be the injective hull of $\top Re_i$. The standard modules are $\Delta(3)=Re_3=\overline{\Delta}(3)$, $\Delta(2)=Re_2/Re_3$ and $\Delta(1)=\top Re_1$. As observed in \cite{CM2021}, $\Hom_R(R/Re_3R, R)$ is the proper standard module $\overline{\Delta}(2)$, which has infinite projective dimension as $R/Re_3R$-module.
	The minimal injective resolution of $Re_3$ is then $I_3^\bullet: 0\rightarrow Re_3\rightarrow I(1)\rightarrow I(1)\rightarrow I(3)\rightarrow 0$ while the minimal injective resolution of $Re_2$ is
	$I_2^\bullet: 0 \rightarrow Re_2\rightarrow Re_1\rightarrow Re_1\rightarrow I(2)\rightarrow 0$.

        Denote the indecomposable projective $R/Re_3R$-modules by $P_{R/Re_3R}(i)$ for $i=1, 2$.
        So $$\rHom_{R}(R/Re_3R, Re_3)\cong \Hom_R(R/Re_3R, I_3^\bullet)$$ is the complex
	$\cdots 0\rightarrow P_{R/Re_3R}(1)\xrightarrow{f} P_{R/Re_3R}(2) \rightarrow 0\rightarrow \cdots$, where the image of $f$ has dimension $3$ and the non-zero modules are placed in degrees $0$ and $1$. Thus, $\rHom_R(R/Re_3R, Re_3)\in \Kb{\pmodcat{R/Re_3R}}$.

        Analogously $\Hom_R(R/Re_3R, I(2))$ is the injective envelope of $\top Re_2$ over $R/Re_3R$ and it has projective dimension two. Therefore, $\rHom_R(R/Re_3R, Re_2)\in  \Kb{\pmodcat{R/Re_3R}}$ and consequently $\rHom_R(R/Re_3R, R)\in \Kb{\pmodcat{R/Re_3R}}$.
	Moreover, $R$ is Iwanaga--Gorenstein and the stratification can be summarised in terms of the following recollements:
\[
\xymatrix@C=1.2cm@R=1.3cm{
	\Kb{\pmodcat{k}} 
	\ar[d] 
	\ar@/^1.5pc/[d] 
	\ar@/_1.5pc/[d] 
	& & \\
	\Kb{\pmodcat{R/Re_3R}} 
	\ar[rr] 
	\ar@/^1.5pc/[rr] 
	\ar@/_1.5pc/[rr] 
	\ar[d] 
	\ar@/^1.5pc/[d] 
	\ar@/_1.5pc/[d] 
	&& \Kb{\pmodcat{R}} 
	\ar@/^1.5pc/[ll] 
	\ar@/_1.5pc/[ll] 
	\ar[rr] 
	\ar@/^1.5pc/[rr] 
	\ar@/_1.5pc/[rr] 
	&& \Kb{\pmodcat{k}} 
	\ar@/^1.5pc/[ll] 
	\ar@/_1.5pc/[ll] \\
	\Kb{\pmodcat{k[x]/(x^2)}} 
	\ar@/^1.5pc/[u] 
	\ar@/_1.5pc/[u] 
	& & & &
}
\]
	
\end{Example}

\medskip

\medskip

\medskip

{\bf AI Declaration.} No generative AI tools were used in this research or
the preparation of this manuscript.

\bigskip

{\bf Acknowledgements.}
An earlier version of this work was presented in the Bielefeld BIREP
seminar on 17.04.2026 and in the Porto Algebra, Combinatorics and Number
Theory seminar on 03.03.2026. The second-named author thanks for the
feedback received during both events.

The first-named author would like to express her
gratitude to the second-named and the third-named
author for their hospitality during her visit to the University of
Stuttgart.

\bigskip

\bibliographystyle{aomalpha}
\def\doi#1.{\unskip\ignorespaces}
\footnotesize{\bibliography{refData}}

\bigskip
\bigskip

\footnotesize{

\noindent
{\sc Yiping Chen, School of Mathematics and Statistics, Wuhan University,
Wuhan, Hubei, 430072,  China
\\
{\tt Email: ypchen@whu.edu.cn}
\medskip

\noindent
{\sc Tiago Cruz, Institute of Algebra and Number Theory, University of
Stuttgart, Pfaffenwaldring 57, 70569 Stuttgart, Germany}
\\
{\tt Email: tiago.cruz@mathematik.uni-stuttgart.de}
\medskip

\noindent
{\sc Steffen Koenig, Institute of Algebra and Number Theory, University of
Stuttgart, Pfaffenwaldring 57, 70569 Stuttgart, Germany}
\\
{\tt Email: skoenig@mathematik.uni-stuttgart.de}
\medskip

}

\end{document}